\documentclass[a4paper,12pt]{amsart}
\usepackage{amssymb,latexsym,graphicx,amsxtra,amsmath,amsthm}
\usepackage{hyperref}

\newcommand{\fr}{\mathcal F}
\newcommand{\ity}{\infty}
\newcommand{\C}{\mathbb{C}}

\numberwithin{equation}{section}
\newtheorem{theorem}{Theorem}[section]
\newtheorem{lemma}[theorem]{Lemma}

\theoremstyle{remark}

\newtheorem{example}[theorem]{Example}
\newtheorem{definition}[theorem]{Definition}
\newtheoremstyle{defi}{10pt}{10pt}{\rm}{\parindent}{\bf}{. }{ }{}

\thanks {The research work of the first author is supported by junior research fellowship from UGC India. }

\begin{document}
\title[NORMALITY AND SHARING FUNCTIONS ]{NORMALITY AND SHARING FUNCTIONS}
\author[G.Datt]{Gopal Datt}
\address{Department of Mathematics, University of Delhi,
Delhi--110 007, India}
\email{ggopal.datt@gmail.com,datt.gopal@ymail.com }
\author[S. Kumar]{Sanjay Kumar}
\address{Department of Mathematics, Deen Dayal Upadhyaya College, University of Delhi,
Delhi--110 015, India }
\email{sanjpant@gmail.com}
\begin{abstract}{In this article, we prove a normality criterion for a family of meromorphic functions which involves sharing of holomorphic functions. Our result generalizes some of the results of H. H. Chen, M. L. Fang ~\cite{CF 1} and M. Han, Y. Gu ~\cite{HG}.
}
\end{abstract}
\maketitle
\section{Introduction and Main Result}\label{sec1}
We denote the  complex plane by $\C$, and the unit disk $\{z\in \C: \ |z|<1 \}$ by $\Delta$.
\begin{definition}\{~\cite{Schiff}, P.33, 71; ~\cite{Ahl}, P.220, 225\}
   Let $\mathcal D$ be a domain on $ \mathbb{C} $,  and $\mathcal{F}$ be a family of meromorphic functions defined on $\mathcal D$. The family $\mathcal{F}$ is said to be normal in $\mathcal{D}$, if every sequence $\{f_n\}\subset \mathcal{F }$ has a subsequence $\{f_{n_j}\}$ which converges spherically  uniformly on compact subsets of  $\mathcal{D}$, to a meromorphic function or $\infty$.  \\
 \end{definition}

According to Bloch's principle every condition that reduces a meromorphic function in the plane to a constant, makes the family of meromorphic functions in a domain $\mathcal{D}$ normal. Rubel gave four counter examples to Bloch principle. \\

Let $f$ and $g$ be  meromorphic functions in a domain $D$ and  $a\in\mathbb{C}.$ Let  zeros of $f-a$ are zeros of  $g-a$ (ignoring multiplicity), we write $f=a\Rightarrow g=a.$ Hence $f=a\iff g=a$ means that $f-a$ and $g-a$ have the same zeros (ignoring multiplicity).
   If $f-a\iff g-a$,   then we say that $f$ and $g$ share the value $z=a$ IM. \{ ~\cite{ccyang}, p. 108\}\\
   
   In this paper, we use the following standard notations of value distribution theory,
\begin{center}
$T(r,f),  m(r,f),  N(r,f),\ldots$.
\end{center}
We denote $S(r,f)$ any function satisfying
\begin{center}
$S(r,f)=o\{T(r,f)\}$,  as $r\rightarrow +\ity,$
\end{center}
 possibly outside of a set with finite measure.\\

In ~\cite{MS}  Mues and Steinmetz proved, {\it{if a non-constant meromorphic  function  $f$ in the plane,  shares three distinct complex numbers $a_1, a_2,  a_3$ with its first order derivative $f'$ then $f \equiv f'$.}}

W. Schwick ~\cite{Sch} was the first who gave a connection between normality and shared values and  proved a theorem related to above result of ~\cite{MS}  which says that  {\it{the family $\fr$ of meromorphic functions on a domain $\mathcal{D}$ is normal if $f$ and $f'$ share $a_1$, $a_2$, $a_3$ for every $f\in \mathcal F$, where $a_1$, $a_2$, $a_3$ are distinct complex numbers.}} Since then many results of normality criteria concerning sharing values have been obtained. In 2001 H. H. Chen and M. L. Fang ~\cite{CF 1} generalized the result of W. Schiwck as follows
\begin{theorem}~\cite{CF 1}
Let $\fr$ be a family of meromorphic functions in a domain $\mathcal{D}$, let $k \ (\geq2)$ be an integer and $a, b, c$ are  complex numbers such that $a\neq b$. If for all $f\in\fr, \ f \ \text{and}\ f^{(k)}$ share $a, b$ and the zeros of $f-c$ are of multiplicity at least $k+1$, then $\fr$  is normal in $\mathcal{D}$.
\end{theorem}

The above result  was further generalized by M. Han and Y. Gu as
\begin{theorem}~\cite{HG}
Let $\fr$ be a family of meromorphic functions in a domain $\mathcal{D}$,  $a, b, c$ are  complex numbers such that $a\neq b$. If for each $f\in\fr$ the zeros of $f-c$ are of multiplicity at least $k+1$, and $f(z)=a$ whenever $f^{(k)}(z)=a$, $f(z)=b$ whenever $f^{(k)}(z)=b$, then $\fr$  is normal in $\mathcal{D}$.
\end{theorem}
One may ask whether one can replace the values $a$ and $b$ by holomorphic functions $a(z)\ \text{and} \ b(z)$ and $f^{(k)}$ by a differential polynomial in $f$.  In this article, we investigate this situation by replacing the values $a$, $b$ by holomorphic functions $a(z)$,  $b(z)$ respectively and $f^{(k)}$ by a linear differential polynomial  $ a_0f^{(k)}(z)+a_1 f(z)$, where  $a_o,  a_1 $ are  complex numbers ,  with $a_0\neq 0$. We define $D(f^{(k)}(z)):=a_0f^{(k)}(z)+a_1f(z)$. \\

Here is our main theorem
 \begin{theorem}\label{1.1}
 Let $\fr$ be a family of meromorphic functions in a domain $\mathcal{D}$, let $k \ (>0)$ be an integer and let $a(z),\ b(z), \ c(z)$ be holomorphic functions such that $a(z)\neq b(z)\ \forall z\in \mathcal{D}$. If for each $f\in \fr$
 \begin{enumerate}
 \item [$(a)$] {the zeros of $f(z)-c(z)$ are of multiplicity $\geq k+1$,  }
 \item [$(b)$] {$f(z)=a(z)$ whenever $D(f^{(k)}(z))=a(z)$},
 \item [$(c)$] {$f(z)=b(z)$ whenever $D(f^{(k)}(z))=b(z)$},

 \end{enumerate}
 then $\fr$ is normal in $\mathcal{D}$.\\
 \end{theorem}

By the following example we observe that the multiplicity restriction on zeros of $f(z)-c(z)$ is sharp in Theorem \ref{1.1}
 \begin{example}
 Let $\mathcal{D} = \Delta =  \{z: | z | < 1\}$, let $k$ be a positive integer, let $a_0=1, a_1=0$, let  $a(z)= z,   b(z) =z+1,  c(z)=0$ (a constant function)  and let
\begin{equation}\notag
\fr=\{nz^k : n=1, 2, 3,\ldots\}.
\end{equation}
Then, for every $f_n(z)\in \fr$, $f_n^{(k)}(z)=nk!$ in $\mathcal{D}$. So $f_n^{(k)}(z)\neq z$ in $\mathcal{D}$ and $f_n^{(k)}(z)\neq z+1$ $\mathcal{D}$. Clearly all conditions of theorem are satisfied but $\fr$ is not normal in $\mathcal{D}$.\\
\end{example}

The following example shows that $a(z_0)$  and $b(z_0)$ can not be equal for some $z_0$ in $\mathcal{D}.$
\begin{example}
Let $k(>0)$ be an integer, let $\mathcal{D} = \Delta =  \{z: | z | < 1\}$, let $a_0=1, a_1=0$, let $a(z)=z, b(z)=2z, c(z)=0 $ (a constant  function). \{So at $z_0=0$ $a(z_0)=b(z_0).$\} and let
 \begin{equation}\notag
 \fr=\{nz^{k+1} : n=1, 2, 3,\ldots\}.
  \end{equation}
  Then, $f_n$ has zero of multiplicity $k+1$ and $f_n^{(k)}(z)=n(k+1)!z$. Clearly all conditions of theorem are satisfied but $\fr$ is not normal in $\mathcal{D}$.\\
\end{example}

The following example shows that the condition that  {$f(z)=a(z)$ whenever $D(f^{(k)}(z))=a(z)$} can not be dropped.
\begin{example}\{~\cite{JFC}, ~\cite{CF}\}
Let $k(>0)$ be an integer, let $\mathcal{D}=\Delta=\{z: |z|<1\}$, $a_0=1, a_1=0$, let $a(z)=0, b(z)=1, c(z)=0$ and let $\fr=\{f_n:n=1, 2, 3, \ldots\},$ where
\begin{equation}\notag
f_n(z)=\frac{(z+\frac{1}{(k+1)n})^{k+1}}{k!(z+\frac{1}{n})}= \frac{z^k}{k!}+p_{k-2}(z)+\frac{(-\frac{k}{n})^{k+1}}{(k+1)^{k+1}k!(z+\frac{1}{n})},
\end{equation}
where $p_{k-2}(z)$ is a polynomial of degree $k-2$. Then, for every $f_n\in\fr,$ $f_n$ has zero of multiplicity $k+1$, and
\begin{equation}\notag
f_n^{(k)}(z)=1-\frac{(\frac{k}{n})^{k+1}}{(k+1)^{k+1}(z\frac{1}{n})^{k+1}},
\end{equation}
Thus $f_n^{(k)}(z)\neq 1$ this shows that $f(z)=b(z)$ whenever $D(f^{(k)}(z))=b(z)$ also $f_n^{(k)}$ has exactly $k+1$ distinct zeros and $f_n$ has only one zero so condition $f(z)=a(z)$ whenever $D(f^{(k)}(z))=a(z)$ is not satisfied. Since $f_n $ takes on the values $0$ and $\ity$ in any fixed neighborhood of $0$, if $n$ is sufficiently large, so $\fr$ fails to be normal in $\mathcal{D}$.
\end{example}
We will use the tools of M. Han and Y. Gu ~\cite{HG} which they used in their paper.

\section{Some Lemmas}
 In order to prove our results we need  the following Lemmas.\\

\begin{lemma}\{~\cite{Zalc}, p. 216; ~\cite{Zalc 1}, p. 814\}\label{lem1}(Zalcman's lemma)\\Let $\mathcal F$ be a family of meromorphic  functions in the unit disk  $\Delta$, with the property that for every function $f\in \mathcal F,$  the zeros of $f$ are of multiplicity at least k. If $\mathcal F$ is not normal at $z_0$ in $\Delta$, then for 0 $\leq \alpha <k$, there exist
\begin{enumerate}
\item{ a sequence of complex numbers $z_n \rightarrow z_0$, $|z_n|<r<1$},
\item{ a sequence of functions $f_n\in \mathcal F$ },
\item{ a sequence of positive numbers $\rho_n \rightarrow 0$},
\end{enumerate}
such that $g_n(\zeta)=\rho_n^{-\alpha}f_n(z_n+\rho_n\zeta) $ converges to a non-constant meromorphic function $g$ on $\C$. Moreover $g$ is of order at most two . \\
\end{lemma}

\begin{lemma}\{~\cite{WF}, p. 21; ~\cite{BE}, p. 360\}\label{lem2}
Let $f(z)$ be a transcendental meromorphic function of finite order, $k$ be a positive integer. If the zeros of $f(z)$ are of multiplicity at least $k+1$, then $f^{(k)}(z)-b$ has infinitely many zeros for any non-zero complex number $b$.\\
\end{lemma}

\begin{lemma}\{~\cite{WF}, p. 22\}\label{lem3}
Let $k (>0)$ be an integer, let $f(z)=a_nz^n+a_{n-1}z^{n-1}+\ldots +a_0+\frac{q(z)}{p(z)},$ where $a_i, i=0, 1, 2, \ldots, n$ are constants with $a_n\neq0, q(z)$ and $p(z)$ are co-prime polynomials with $\text{deg}\ q(z)< \text{deg}\ p(z).$ If $f^{(k)}(z)\neq1,$ then
\begin{equation}\notag
f(z)=\frac{z^k}{k!}+\ldots +a_0+\frac{1}{(a z+b)^n},
\end{equation}
where $a\neq0$. If all zeros of $f(z)$ have multiplicity at least $k+1$, then
\begin{equation}\notag
f(z)=\frac{(cz+d)^{k+1}}{az+b},
\end{equation}
where $c (\neq0), d$ are constants.\\  
\end{lemma}

\begin{lemma}\{~\cite{ccyang}, p. 43; ~\cite{Yang}, p. 110\}\label{lem4}(Hayman's Inequality)
Let $k (>0)$ be an integer, suppose that $f(z)$ is a transcendental meromorphic function in the complex plane. Then
\begin{equation}\notag
T(r,f)< (2+\frac{1}{k})N(r,\frac{1}{f})+(2+\frac{2}{k})N(r,\frac{1}{f^{(k)}-1})+S(r,f).
\end{equation}
\end{lemma}

\section{Proof of Theorem \label{theorem 1.1}}
 \begin{proof}
 Since normality is a local property, we assume that $D=\Delta=\{z:|z|<1\}.$
Suppose  $\fr$ is not normal in $D$. Without loss of generality we assume that $\fr$ is not normal at the point $z_0=0$ in $\Delta$.  Then by Lemma \ref{lem1}, there exist
\begin{enumerate}
\item { a sequence of complex numbers $z_n \rightarrow z_0$, $|z_n|<r<1$},
\item { a sequence of functions $f_n\in \mathcal F$ and}
\item { a sequence of positive numbers $\rho_n \rightarrow 0$,}
\end{enumerate}
such that $g_n(\zeta)=\frac{f_n (z_n+\rho_n\zeta)-c(z_n+\rho_n\zeta)}{\rho_n^k}$ converges locally uniformly to a non-constant meromorphic  function $g(\zeta)$ in $\C$ and  the zeros of
$g(\zeta)$ are of multiplicity at least $k+1$.  Moreover $g$ is of order at most two. \\

\textbf{Case 1} When $c(z_0)=a(z_0)$, we claim that
\begin{enumerate}
\item [(i)] {$g^{(k)}(\zeta)\neq B$, where $B=\frac{b(z_0)-a_0c^{(k)}(z_0)+a_1c(z_0)}{a_0}$ \ is a constant. }
\item [(ii)] {$g^{(k)}(\zeta)=A \Rightarrow g(\zeta)=0$, where $A=\frac{a(z_0)-a_0c^{(k)}(z_0)-a_1(z_0)}{a_0}$ \ is a constant.}
\end{enumerate}
Clearly $g^{(k)}(\zeta)\not\equiv B$ as zeros of $g(\zeta)$ are of multiplicity at least $k+1$. Suppose $g^{(k)}(\zeta_0)=B,$ then by Hurwitz's theorem there exist $\zeta_n;  \  \zeta_n\rightarrow \zeta_0$ such that
\begin{equation}\notag
 g_n^{(k)}(\zeta_n)=f_n^{(k)}(z_n+\rho_n\zeta_n)-c^{(k)}(z_n+\rho_n\zeta_n)=B.
 \end{equation}
 Now, consider
 \begin{align}
 \frac{D(f^{(k)}(z_n+\rho_n\zeta_n)-b(z_n+\rho_n\zeta_n)}{a_0}&=\frac{a_0f_n^{(k)}(z_n+\rho_n\zeta_n)+a_1f_n(z_n+\rho_n\zeta_n)-b(z_n+\rho_n\zeta_n)}{a_0}\notag\\
 &=f_n^{(k)}(z_n+\rho_n\zeta_n)+ \frac{a_1f_n(z_n+\rho_n\zeta_n)-b(z_n+\rho_n\zeta_n)}{a_0}\notag\\
 &=f_n^{(k)}(z_n+\rho_n\zeta_n)+\frac{a_1}{a_0}[\rho_n^kg_n(\zeta_n)+c(z_n+\rho_n\zeta_n)]\notag\\
 &-\frac{b(z_n+\rho_n\zeta_n)}{a_0}\notag\\
 &\rightarrow g^{(k)}(\zeta_0)+\frac{a_1}{a_0}c(z_0)-\frac{b(z_0)}{a_0}+c^{(k)}(z_0)\notag\\
 &=g^{(k)}(\zeta_0)-\frac{b(z_0)-a_0c^{(k)}(z_0)+a_1c(z_0)}{a_0}\label{eq 1}.
 \end{align}
 Now, it follows from \ref{eq 1} and condition $(c)$ of the theorem \ref{1.1} that
 \begin{equation}\notag
 f(z_n+\rho_n\zeta_n)=b(z_n+\rho_n\zeta_n).
 \end{equation}
  Thus,
  \begin{align*}
 g_n(\zeta_n)&=\frac{f_n(z_n+\rho_n\zeta_n)-c(z_n+\rho_n\zeta_n)}{\rho_n^k}\\
 &=\frac{b(z_n+\rho_n\zeta_n)-c(z_n+\rho_n\zeta_n)}{\rho_n^k}\rightarrow \ity.
 \end{align*}
 $i.e. \ g(\zeta_0)=\ity$ which contradicts $g^{(k)}(\zeta_0)=B$. Thus $g^{(k)}(\zeta)\neq B.$\\

 Now, we prove claim (ii)
 Suppose $g^{(k)}(\zeta_0)=A$, by Hurwitz's theorem, there exist exist $\zeta_n;  \  \zeta_n\rightarrow \zeta_0$ such that
\begin{equation}\notag
 g_n^{(k)}(\zeta_n)=f_n^{(k)}(z_n+\rho_n\zeta_n)-c^{(k)}(z_n+\rho_n\zeta_n)=A.\\
 \end{equation}

 Now, consider
 \begin{align}
 \frac{D(f^{(k)}(z_n+\rho_n\zeta_n)-a(z_n+\rho_n\zeta_n)}{a_0}&=\frac{a_0f_n^{(k)}(z_n+\rho_n\zeta_n)+a_1f_n(z_n+\rho_n\zeta_n)-a(z_n+\rho_n\zeta_n)}{a_0}\notag\\
 &\rightarrow g^{(k)}(\zeta_0)+\frac{a_1}{a_0}c(z_0)-\frac{a(z_0)}{a_0}+c^{(k)}(z_0)\notag\\
 &=g^{(k)}(\zeta_0)-\frac{a(z_0)-a_0c^{(k)}(z_0)+a_1c(z_0)}{a_0}\label{eq 2}.
 \end{align}
 So, it follows from \eqref{eq 2} and condition $(b)$ of the theorem \ref{1.1} that
  \begin{equation}\notag
 f(z_n+\rho_n\zeta_n)=a(z_n+\rho_n\zeta_n).
 \end{equation}
  Thus,
  \begin{align*}
 g_n(\zeta_n)&=\frac{f_n(z_n+\rho_n\zeta_n)-c(z_n+\rho_n\zeta_n)}{\rho_n^k}\\
 &=\frac{a(z_n+\rho_n\zeta_n)-c(z_n+\rho_n\zeta_n)}{\rho_n^k}=0.
 \end{align*}
So $g(\zeta_0)=0.$\\

 \textbf{Subcase 1.1} If $B\neq0$, it is evident from Lemma \ref{lem2}, that $g(\zeta)$  is not a transcendental meromorphic function as $g^{(k)}(\zeta)\neq B.$ Hence, $g(\zeta)$ is a non-constant rational function. Then, by Lemma \ref{lem3}, we have
 \begin{equation}\label{eq 3}
 g(\zeta)=\frac{B\zeta^k}{k!}+\ldots+b_0+\frac{1}{(\alpha\zeta+\beta)^n}.
 \end{equation}
 Thus
 \begin{equation}\label{eq 4}
 g^{(k)}(\zeta)=B+\frac{D}{(\alpha\zeta+\beta)^{n+k}},
 \end{equation}
where $D$ is a constant, by \eqref{eq 3} and \eqref{eq 4} we conclude that the number of the zeros of $g(\zeta)$ is $k+n$, and by, Hurwitz's theorem multiplicity of zeros are at least $k+1$, so the number of the distinct zeros of $g(\zeta)$ is   at most $\frac{k+n}{k+1}.$ It is simple to check that the zeros of $g^{(k)}(\zeta)-A$  are of multiplicity $1,$ so the number of the distinct zeros of $g^{(k)}$ is $k+n$, which does not hold claim (ii).\\

\textbf{Subcase 1.2} If $B=0$, then by, claim (i) $g^{(k)}(\zeta)\neq 0$. This shows that either $g(\zeta)\neq 0$ or the zeros if $g(\zeta) $ are of multiplicity  at most $ k$, which contradicts that the zeros of $g(\zeta)$ are of multiplicity at least $k+1$. So we have $g(\zeta)\neq0$. By claim (ii) we get $g^{(k)}(\zeta)\neq A,$ where $A\neq0$. Since $g(\zeta)\neq0, \ g^{(k)}(\zeta)\neq A, \ A\neq0$, this implies $N(r,\frac{1}{g})=0 \ \text{and}\ N(r,\frac{1}{g^{(k)}-A})=0$, so by Lemma \ref{lem4} (Hayman's inequality) we get that $g(\zeta)$ is a constant, which is a contradiction.\\

\textbf{Case 2} When $c(z_0)\neq a(z_0)$ and $c(z_0)\neq b(z_0)$. By Lemma \ref{lem1}, there exist
\begin{enumerate}
\item { a sequence of complex numbers $z_n \rightarrow z_0$, $|z_n|<r<1$},
\item { a sequence of functions $f_n\in \mathcal F$ and}
\item { a sequence of positive numbers $\rho_n \rightarrow 0$,}
\end{enumerate}
such that $g_n(\zeta)=\frac{f_n (z_n+\rho_n\zeta)-c(z_n+\rho_n\zeta)}{\rho_n^k}$ converges locally uniformly to a non-constant meromorphic  function $g(\zeta)$ in $\C$  The zeros of
$g(\zeta)$ are of multiplicity at least $k+1$.  Moreover $g$ is of order at most two. \\
we claim that
\begin{enumerate}
\item [(a)] {$g^{(k)}(\zeta)\neq B$, where $B=\frac{b(z_0)-a_0c^{(k)}(z_0)+a_1c(z_0)}{a_0}$ \ is a constant. }
\item [(b)] {$g^{(k)}(\zeta)\neq A$, where $A=\frac{a(z_0)-a_0c^{(k)}(z_0)-a_1(z_0)}{a_0}$ \ is a constant.}
\end{enumerate}
Using the method of Case $1$  these claims can be proven. By Lemma \ref{lem2}, we have $g(\zeta)$ is a non-constant rational function, so  is $g^{(k)}(\zeta)$ and by claim $g^{(k)}(\zeta)$ has two omitted values, which is a contradiction. \\

\textbf{Case 3} when $c(z_0)=b(z_0)$. this case is similar to case 1.

 \end{proof}


\begin{thebibliography}{00}

\bibitem{Ahl} L. V. Ahlfors, Complex Analysis, Third edition, McGraw-Hill, 1979.
\bibitem{BE} W. Bergweiler, A. Eremenko, On The Singularities of The Inverse to a Meromorphic Function of Finite Order, Rev. Mat. Iberoamericana, \textbf{11} (1995), No. 2, 355--372.
\bibitem{CF 1} H. H. Chen, M. L. Fang, Shared Values and Normal Families of Meromorphic Functions, J.  Math. Anal. Appl., \textbf{260} (2001), 124--132. 
\bibitem{JFC}J. F. Chen, Normal Families and Shared Sets of Meromorphic Functions, Rocky Mountain Journal of Mathematics, \textbf{41} (2011), No. 1, 37--43.
\bibitem{CF} J. F. Chen, M. L. Fang, Normal Families And Shared Functions of Meromorphic Functions, Israel Journal of Mathematics  \textbf{180} (2010), 129--142.
\bibitem{MLF} M. L. Fang, Picard Values and Normality criterion, Bull. Korean Math. Soc. \textbf{38} (2001), No. 2, 379--387.
\bibitem{HG} M. Han, Y. Gu, The Normal Family of Meromophic Functions, Acta Mathematica Scientia, \textbf{28 B} (2008), No. 4, 759--762.
\bibitem{Hay}W. K. Hayman, Meromorphic Functions, Clarendon Press, Oxford, 1964.
\bibitem{MS}E. Mues and N. Steinmentz, Meromorphe Functionen, die mit ihrer Ableitung Werte teilen, Manuscripta Mathematica \textbf{29} (1979), 195--206.
\bibitem{XPLZ} X. Pang and L. Zalcman, Normality and Shared values, Ark. Mat. \textbf{38} (2000), 171--182.
\bibitem{Schiff} J. Schiff,  Normal Families, Springer-Verlag, Berlin, 1993.
\bibitem{Sch}W. Schwick, Sharing Values and Normality, Archiv der Mathematik \textbf{59},(1992), 50--54.
\bibitem{WF}Y. F. Wang, M. L. Fang, Picard Values and Normal Families of Meromorphic Functions With Multiple Zeros, Acta Math. Sin., new series, \textbf{14}  (1998), No. 1, 17--26.
\bibitem{XX}J. Xia, Y. Xu, Normality Criterion Concerning Sharing Functions II, Bull. Malays. Math. Sci. Soc.(2) \textbf{33} (2010), no. 3, 479--486.
\bibitem{YX}Y. Xu, Normality Criterion concerning Sharing Functions, Houstan Journal of Mathematics \textbf{32} (2006),  No. 3,  945--954.
\bibitem{ccyang} C. C. Yang, H. X. Yi, Uniqueness theory of meromorphic functions, Science Press/ Kluwer Academic Publishers, 2003.
\bibitem{Yang}L. Yang, Value Distribution Theory, Springer- Verlag,Berlin, 1993.
\bibitem{Zalc 1}L. Zalcman, A heuristic principle in complex function theory, The Amer. Math. Monthly \textbf{82} (1975), 813--817.
\bibitem{Zalc}L. Zalcman, Normal Families: New perspective, Bulletin of American Mathematical Society \textbf{35} (1998), 215--230.







\end{thebibliography}
\end{document}